\documentclass[a4paper,10pt]{article}

\usepackage{amssymb, amsmath, amsthm, latexsym, verbatim, enumerate,array}

\renewcommand \a{\alpha}

\newcommand \K{\delta}

\newcommand \la{\lambda}
\newcommand \id{\mathrm{id}}
\newcommand \Rn{\mathbb R^n}
\newcommand \RN{\mathbb R^N}
\newcommand \rk{\mathrm{rk}}

\newcommand \Span{\mathrm{Span}}
\newcommand \Tr{\mathrm{Tr} \, }

\newcommand \diag{\mathrm{diag}}
\newcommand \mm{\mathbf{m}}

\newcommand \<{\langle}
\renewcommand \>{\rangle}

\newcommand \g{\mathfrak{g}}
\newcommand \p{\mathfrak{p}}

\newcommand \n{\mathfrak{n}}
\newcommand \m{\mathfrak{m}}
\newcommand \z{\mathfrak{z}}
\renewcommand \k{\mathfrak{k}}
\renewcommand \b{\mathfrak{b}}
\newcommand \ag{\mathfrak{a}}
\newcommand \ad{\mathrm{ad}}
\newcommand \ric{\mathrm{ric}}
\newcommand \Ric{\mathrm{Ric}}
\newcommand \Der{\mathrm{Der}}
\newcommand \End{\mathrm{End}}

\newcommand \GL{\mathrm{GL}}
\newcommand \Conv{\mathrm{Conv}}

\theoremstyle{plane}
\newtheorem{theorem}{Theorem}
\newtheorem*{theorem*}{Theorem}
\newtheorem*{corollary*}{Corollary}
\newtheorem{lemma}{Lemma}

\newtheorem*{proposition*}{Proposition}

\newtheorem*{namedtheorem}{\theoremname}
\newcommand{\theoremname}{te}

\theoremstyle{remark}
\newtheorem{remark}{Remark}
\newtheorem{example}{Example}

\begin{document}

\title{Einstein solvmanifolds with a simple Einstein derivation}

\author{Y.Nikolayevsky}

%\date{\today}
%\keywords
%\subjclass Primary: 53C25?

\maketitle

\begin{abstract}
The structure of a solvable Lie groups admitting an Einstein left-invariant metric is, in a sense, completely
determined by the nilradical of its Lie algebra. We give an easy-to-check necessary and sufficient condition for a
nilpotent algebra to be an Einstein nilradical whose Einstein derivation has simple eigenvalues. As an application,
we classify filiform Einstein nilradicals (modulo known classification results on filiform graded Lie algebras).
\end{abstract}

\section{Introduction}
\label{s:intro}

In this paper, we study Riemannian homogeneous spaces with an Einstein metric of negative scalar curvature.
The major open conjecture here is the \emph{Alekseevski Conjecture} \cite{A1} asserting that an Einstein homogeneous
Riemannian space
with negative scalar curvature admits a simply transitive solvable isometry group. This is equivalent to saying
that any such space is a \emph{solvmanifold}, a solvable Lie group with a left-invariant Riemannian metric satisfying
the Einstein condition.

Even assuming the Alekseevski Conjecture, the classification (or even the description) of Einstein solvmanifolds is
quite a complicated problem. Very recently, an amazing progress in this direction was achieved by J.Lauret \cite{L4}
who proved that any Einstein solvmanifold is \emph{standard}.
This means that the metric solvable Lie algebra $\g$ of such a solvmanifold has the following property: the orthogonal
complement $\ag$ to the derived algebra of $\g$ is abelian. The systematic study of standard Einstein solvmanifolds
(and the term ``standard") originated from the paper \cite{H}. In particular, it is shown there that
any non-unimodular Einstein solvmanifold admits a rank-one reduction (any unimodular Einstein solvmanifold is flat by
\cite{DM}). On the Lie algebra
level, this means that if $\g$ is an Einstein metric solvable Lie algebra and $\g=\ag \oplus \n$
(orthogonal sum), with $\n$ the nilradical of $\g$, then there exists a one-dimensional subspace $\ag_1 \subset \ag$ such
that $\g_1 = \ag_1 \oplus \n$ with the induced inner product is again an Einstein metric solvable Lie algebra.
What is more, $\g_1$ is essentially uniquely determined by $\n$, and all the Einstein metric solvable Lie
algebras with the nilradical $\n$ can be obtained from $\g_1$ via a known procedure (by adjoining appropriate
derivations).

In particular, the geometry (and the algebra) of an Einstein metric solvable Lie algebra is completely encoded
in its nilradical. A nilpotent Lie algebra which can be a nilradical of a
an Einstein metric solvable Lie algebra is called an \emph{Einstein nilradical}.

A vector $H$ spanning $\ag_1$ and scaled in such a way that $\|H\|^2 = \Tr \ad_H$ is called
\emph{the mean curvature vector},
and the restriction of the derivation $\ad_H$ to $\n$ the \emph{Einstein derivation}. As it is proved in \cite{H},
the Einstein derivation is always semisimple, and, up to scaling, its eigenvalues are natural numbers.
The
ordered set of eigenvalues $\la_i \in \mathbb{N}$ of the (appropriately scaled) Einstein derivation, together with the
multiplicities $d_i$, is called the \emph{eigenvalue type} of an Einstein solvmanifold (written as
$(\la_1 < \ldots < \la_p; \, d_1, \ldots, d_p)$).

Considerable efforts were made towards the classification of Einstein nilradicals (Einstein solvmanifolds) of a given
eigenvalue type. While in the trivial case of eigenvalue type $(1 ; n)$ there is only one possible Einstein nilradical
(the abelian one; the corresponding solvmanifold is the hyperbolic space), even the eigenvalue type $(1 < 2 ; d_1, d_2)$
is far from being completely understood (see \cite{GK, H}).

In this paper, we study solvmanifolds with the simple eigenvalue type, $(\la_1 < \ldots < \la_n; \, 1, \ldots, 1)$,
and give an easy-to-check necessary and sufficient condition to determine, whether a given nilpotent Lie algebra is
an Einstein nilradical with such an eigenvalue type. We have to answer two questions: %Our task here is two-fold:
\begin{enumerate}[(i)]
  \item \label{it:1}
  given an arbitrary nilpotent Lie algebra, how to recognize the eigenvalue type of its Einstein solvable
  extension (if for some inner product such an extension exists); in particular, how to determine, whether such an
  extension, if it exists, has a simple eigenvalue type;
  \item \label{it:2}
  if, according to (\ref{it:1}), a nilpotent Lie algebra may potentially admit an Einstein metric solvable extension
  with a simple eigenvalue type, whether it actually does.
\end{enumerate}

The answer to (\ref{it:1}) is given by constructing the
pre-Einstein derivation introduced in \cite{N2}. A semisimple derivation $\phi$ of a nilpotent Lie algebra $\n$ is
called \emph{pre-Einstein}, if
$\Tr \, (\phi \circ \psi) = \Tr \, \psi$, for any $\psi \in \Der (\n)$. A pre-Einstein derivation always exists and
is unique up to conjugation, and its eigenvalues are rational. As it is shown in \cite{N2}, the condition
``$\phi > 0$ and $\ad_{\phi} \ge 0$" is necessary (but not sufficient) for $\n$ to be an Einstein nilradical
(``$A > 0$" means that all the eigenvalues of the operator $A$ are real and positive).
What is more, if $\n$ \emph{is} an Einstein nilradical, then its Einstein derivation is positively proportional to
$\phi$ and the eigenvalue type of the corresponding Einstein metric solvable Lie algebra is just the
set of the eigenvalues of $\phi$ and their multiplicities (see Section~\ref{s:pre} for details).

Our main result of is Theorem~\ref{t:one} below, answering (\ref{it:2}).
Let $\n$ be a nilpotent Lie algebra of dimension $n$, with $\phi$ a pre-Einstein derivation.
Suppose that all the eigenvalues of $\phi$ are simple. Let $e_i$ be the basis of eigenvectors for $\phi$ and let
$[e_i, e_j] = \sum_{k=1}^n C_{ij}^k e_k$ (note that for every pair $(i, j)$, no more than one of the $C_{ij}^k$ is
nonzero). In a Euclidean space $\Rn$ with the inner product $(\cdot, \cdot)$ and an orthonormal basis
$f_1, \ldots, f_n$, define the finite subset $\mathbf{F}=\{\a_{ij}^k= f_i+f_j-f_k: C_{ij}^k \ne 0\}$. Let $L$ be the
affine span of $\mathbf{F}$, the smallest affine subspace of $\Rn$ containing $\mathbf{F}$.

\begin{theorem}\label{t:one}
Let $\n$ be a nilpotent Lie algebra whose pre-Einstein derivation has all the eigenvalues simple.
The algebra $\n$ is an Einstein nilradical if and only if the projection of the origin of $\Rn$ to $L$ lies in the
interior of the convex hull of $\mathbf{F}$.
\end{theorem}

\begin{remark}\label{rem:payne}
Denote $N = \# \mathbf{F}$ and
introduce an $n \times N$ matrix $Y$ whose vector-columns are the vectors $\a_{ij}^k$ in some fixed order. Define
the vector $[1]_N =(1,\ldots, 1)^t \in \mathbb{R}^N$ and an $N \times N$ matrix $U=Y^tY$. One can rephrase
Theorem~\ref{t:one} as follows: \emph{a nilpotent Lie algebra $\n$ whose pre-Einstein derivation has all its eigenvalues
simple is an Einstein nilradical if and only if there exists a vector $v \in \mathbb{R}^N$ all of whose coordinates
are positive such that}
\begin{equation}\label{eq:ut=1}
Uv=[1]_N.
\end{equation}
By the result of \cite[Theorem 1]{P}, a metric nilpotent Lie algebra is nilsoliton, if and only if
equation~\eqref{eq:ut=1} holds with respect to the basis of Ricci eigenvectors.
Note that there are two fundamental differences between Theorem~\ref{t:one} and this result. First of all, our $\n$ is
just a Lie algebra, no inner product is present. Secondly, and more importantly, the set $S$ of vectors $v$ satisfying
\eqref{eq:ut=1} is the interior of a convex polyhedron in an affine subspace of $\RN$ (if it is nonempty). If
$v=(v_{ij}^k)$ is a point from $S$, the skew-symmetric bilinear map on the linear space of $\n$ defined by
$[e_i, e_j] = \pm (v_{ij}^k)^{1/2} e_k$ (where $i < j$ and $k$ is defined by the condition
that $C_{ij}^k \ne 0$) may not even be a Lie bracket, and if it is, there is no apparent reason why the resulting Lie
algebra has to be isomorphic to $\n$ (say, when $S$ has a positive dimension). The fact that the point $v$ can be chosen
``in the correct way" is the main content of Theorem~\ref{t:one}. It should be noted, however, that one can hardly expect
to find the nilsoliton inner product explicitly, except in the cases when $\n$ is particularly nice.
\end{remark}

\bigskip

The class of nilpotent algebras whose pre-Einstein derivation has all its eigenvalues simple is rather broad (see the
end of Section~\ref{s:proof} for some examples). As an application
of Theorem~\ref{t:one}, we give the classification of filiform Einstein nilradicals (a nilpotent Lie
algebra $\n$ of dimension $n$ is called a \emph{filiform}, if its descending central series has the maximal possible
length, $n-1$).

As any filiform algebra $\n$ is generated by two elements, its \emph{rank}, the dimension of the maximal torus of
derivations (the maximal abelian subalgebra of $\Der (\n)$ consisting of semisimple elements), is at most two.
Most of filiform algebras of dimension $n \ge 8$ are characteristically nilpotent \cite{GH1}, that is, have
rank zero. Such algebras do not admit any gradation at all, hence cannot be Einstein nilradicals.
There are two series of filiform algebras of rank two: $\m_0(n)$ (given by the relations
$[e_1, e_i] = e_{i+1}$, $i=2, \ldots n-1$) and
$\m_1(n), \; n$ is even (the relations of $\m_0(n)$ and $[e_i,e_{n-i+1}] = (-1)^i e_n, \; i = 2, \ldots , n - 1$)
\cite{V, GH1}. Both of them are Einstein nilradicals: for $\m_0(n)$, this is proved in \cite[Theorem 4.2]{L2}
(see also \cite[Theorem 27]{P}), for $\m_1(n)$ in \cite[Theorem 37]{P}.

Less is known about filiform algebras of rank one. According to \cite[Th\'{e}or\`{e}me 2]{GH1}, there are two
classes of such algebras admitting a positive gradation: $A_r$ and $B_r$, with $r \ge 2$.
An $n$-dimensional algebra $\n$ of the class $A_r$ ($n \ge r+3$) is given by the relations of $\m_0(n)$ and
$[e_i, e_j] = c_{ij} e_{i+j+r-2}$ for $i,j \ge 2, \; i+j \le n+2-r$.
An $n$-dimensional algebra $\n$ of the class $B_r$ ($n \ge r+3$, $n$ is even)
is given by the relations
$[e_1, e_i] = e_{i+1}, \quad i=2, \ldots n-2, \quad [e_i, e_j] = c_{ij} e_{i+j+r-2} , \quad i,j \ge 2, \; i+j \le n+2-r$.

The complete classification is known only for the algebras of the class $A_2$ (\cite{M}, and independently
in \cite{CJ}): the class $A_2$ consists of five infinite series
and of five one-parameter families $\g_\a(n)$ in dimensions $n=7, \ldots, 11$
(to the best of the author's knowledge, no classification results for $r \ge 3$ appeared in the literature).
Based on that classification, we classify the algebras of the class $B_2$ and prove the following theorem
(here $\mathcal{V}(n)$ is the $n$-dimensional truncated Witt algebra, the others are the members of finite families
listed in the tables in Section~\ref{s:fili}):

\begin{theorem} \label{t:fili}
{\ }

\emph{1.} A filiform algebra $\n \in A_2$ of dimension $n \ge 8$ is an Einstein nilradical if and only if it is
isomorphic either to
$\mathcal{V}(n)$ or to one of the following algebras from Table~\ref{tablega}:
$\g_\a(8), \, \a \ne -2,\quad \g_\a(9), \, \a \ne -2$, $\g_\a(10), \, \a \ne -2,-1,\frac12, \quad \g_\a(11)$.

\emph{2.} The class $B_2$ consists of six algebras $\b(6), \b(9), \b_1(10), \b_2(10), \b_{\pm}(12)$ listed in
Table~\ref{tableb2n}.
All of them are Einstein nilradicals.
\end{theorem}

In assertion 1 of Theorem~\ref{t:fili} we consider only the case $\dim \n \ge 8$, as every nilpotent algebra of
dimension six or lower is an Einstein nilradical (\cite[Theorem 3.1]{W} and \cite[Theorem 5.1]{L2}), and all
seven-dimensional filiform algebras, except for $\m_{0,1}(7) \cong \g_{-2}(7), \, \g_{-1}(7)$, and $\m_2(7)$, are
Einstein  nilradicals by \cite[Theorem 4.1]{LW}.

\bigskip

The paper is organized as follows. In Section~\ref{s:pre} we give the background on Einstein solvmanifolds
and on the momentum map. The proof of Theorem~\ref{t:one} is given in Section~\ref{s:proof}. In Section~\ref{s:fili},
we consider graded filiform algebras and prove Theorem~\ref{t:fili} (in Section~\ref{ss:a2n} for the algebras of
the class $A_2$, and in Section~\ref{ss:b2n} for the algebras of the class $B_2$).

\section{Preliminaries}
\label{s:pre}

For an inner product $\< \cdot, \cdot \>$ on a Lie algebra $\g$, define the \emph{mean curvature vector} $H$ by
$\<H, X\> = \Tr \ad_X$ (clearly, $H$ is orthogonal to the derived algebra of $\g$). For $A \in \End(\g)$, let $A^*$
be its metric adjoint and $S(A) =\frac12 (A +A^*)$ be the symmetric part of $A$. Let $\Ric$ be the Ricci
$(0, 2)$-tensor (a quadratic form) of $(\g, \< \cdot, \cdot \>)$, and $\ric$ be the \emph{Ricci operator}, the symmetric
operator associated to $\Ric$.

The Ricci operator of $(\g, \< \cdot, \cdot \>)$ is implicitly defined by
\begin{equation}\label{eq:riccidef}
\Tr \Bigl(\ric + S(\ad_H) + \frac12 B \Bigr)
\circ A = \frac14 \sum\nolimits_{i,j} \<A[E_i, E_j] - [AE_i, E_j] - [E_i, AE_j], [E_i, E_j]\>,
\end{equation}
for any $A \in \End(\g)$, where $\{E_i\}$ is an orthonormal basis for $\g$, and $B$ is the symmetric operator associated
to the Killing form of $\g$.

If $(\n, \< \cdot, \cdot \>)$ is a nilpotent metric Lie algebra, then $H = 0$ and $B = 0$, so \eqref{eq:riccidef} gives
\begin{equation}\label{eq:riccinil}
\Tr (\ric_{\n} \circ A) = \frac14 \sum\nolimits_{i,j} \<A[E_i, E_j] - [AE_i, E_j] - [E_i, AE_j], [E_i, E_j]\>.
\end{equation}

An inner product on a solvable Lie algebra $\g$ is called \emph{standard}, if the orthogonal complement to
the derived algebra $[\g, \g]$ is abelian. A metric solvable Lie algebra $(\g, \<\cdot,\cdot\>)$ is called
\emph{standard}, if the inner product $\<\cdot,\cdot\>$ is standard.

By the result of \cite{L4}, any Einstein metric solvable Lie algebra must be standard.

As it is proved in \cite{AK}, any Ricci-flat metric solvable Lie algebra is flat. By the result of \cite{DM},
any Einstein metric solvable unimodular Lie algebra is also flat. In what follows, we always assume $\g$
to be nonunimodular ($H \ne 0$), with an inner product of a strictly negative scalar curvature $c \dim \g$.

Any Einstein metric solvable Lie algebra admits a rank-one reduction \cite[Theorem 4.18]{H}. This means that if
$(\g, \< \cdot, \cdot\>)$ is such an algebra, with the nilradical $\n$ and the mean curvature vector $H$, then the
subalgebra $\g_1 = \mathbb{R}H \oplus \n$, with the induced inner product, is also Einstein. What is
more, the derivation $\phi=\ad_{H|\n}:\n \to \n$ is symmetric with respect to the inner product, and all its
eigenvalues belong to $\a \mathbb{N}$ for some constant $\a > 0$. This implies, in particular, that the nilradical $\n$
of an Einstein metric solvable Lie algebra admits an $\mathbb{N}$-gradation defined by the eigenspaces of $\phi$.
As it is proved in \cite[Theorem 3.7]{L1}, a necessary and sufficient condition for a metric nilpotent algebra
$(\n, \< \cdot, \cdot\>)$ to be the nilradical of an Einstein metric solvable Lie algebra is
\begin{equation}\label{eq:ricn}
    \ric_\n = c \, \id_\n + \phi,
\end{equation}
where $c \dim \g < 0$ is the scalar curvature of $(\g, \< \cdot, \cdot\>)$. This equation, in fact, defines
$(\g, \< \cdot, \cdot\>)$ in the following sense: given a metric nilpotent Lie algebra whose Ricci operator
satisfies \eqref{eq:ricn}, with some constant $c < 0$ and some $\phi \in \Der(\n)$, one can define $\g$ as a
one-dimensional extension of $\n$ by $\phi$. For such an extension $\g = \mathbb{R}H \oplus \n, \; \ad_{H|\n} = \phi$,
and the inner product
defined by $\<H, \n \> = 0,\; \|H\|^2 = \Tr \phi$ (and coinciding with the existing one on $\n$) is Einstein, with the
scalar curvature $c \dim \g$. A nilpotent Lie algebra $\n$ which admits an inner product
$\< \cdot, \cdot\>$ and a derivation $\phi$ satisfying \eqref{eq:ricn} is called an \emph{Einstein nilradical}, the
corresponding derivation $\phi$ is called an \emph{Einstein derivation}, and the inner product $\< \cdot, \cdot\>$
the \emph{nilsoliton metric}.

As it is proved in \cite[Theorem 3.5]{L1}, a nilpotent Lie algebra admits no more than one nilsoliton metric, up to
a conjugation and scaling (and hence, an Einstein derivation, if it exists, is unique, up to a conjugation and scaling).
If $\la_1 < \ldots < \la_p$ are the eigenvalues of $\phi$, with $d_1, \ldots, d_p$ the corresponding multiplicities, we
call $(\la_1 < \ldots < \la_p; \; d_1, \ldots, d_p)$
the \emph{eigenvalue type} of the Einstein metric solvable Lie algebra $(\g, \< \cdot, \cdot\>)$. With some abuse of
language, we will also call $(\la_1 < \ldots < \la_p; \; d_1, \ldots, d_p)$ the
\emph{eigenvalue type} of $\phi$.

The main tool in the proof of Theorem~\ref{t:one} is the moment map. Let $G$ be a reductive Lie group, with
$K \subset G$ a maximal compact subgroup. Let $\g= \k \oplus \p$ be the Cartan decomposition of the Lie algebra of $G$,
with $\k$ the Lie algebra of $K$. Suppose $G$ acts on a linear space $V$ endowed with a $K$-invariant inner product
$\<\cdot,\cdot\>$. The action of $G$ is then descends to the projective space $\mathbb{P}V$.

The \emph{moment map} $\mm$ of the action of $G$ on $\mathbb{P}V$ is defined by
\begin{equation}\label{eq:defmoment}
    \mm: \mathbb{P}V \to \p^*, \quad \mm(x)(X)=\frac{1}{\|v\|^2} \frac{d}{dt}_{|t=0}\<\exp(tX).v, v\>,
    \quad \text{for $x =[v], \, X \in \p$}.
\end{equation}

The fact that the moment map can be used to study the nilsoliton metrics was first observed in \cite{L1}, where the
following construction was given.
Let $\mathcal{L}=(\Rn, \<\cdot, \cdot\>)$ be a linear space with the inner product, and let
$V=\mathrm{Hom}(\Lambda^2 \mathcal{L},\mathcal{L})$ be the space of skew-symmetric bilinear maps from $S$ to itself.
Denote $\mathcal{N} \subset V$ the (real algebraic) subset of those $\mu \in V$, which are nilpotent Lie brackets.
The inner product on $V$ is defined in an obvious way: for $\mu_1, \mu_2 \in V$,
$\<\mu_1,\mu_2\>=\sum_{ij}\<\mu_1(e_i,e_j),\mu_2(e_i,e_j)\>$, where $\{e_i\}$ is an orthonormal basis for $\mathcal{L}$.
The group $G=\GL(n)$ acts on $V$ as follows: for $\mu \in V$ and $g \in G$, $g.\mu(X,Y) =g \mu(g^{-1}X, g^{-1}Y)$,
where $X, Y \in \mathcal{L}$ (clearly, $\mathcal{N} \subset V$ is $G$-invariant).
Take $\mathfrak{gl}(n)=\mathfrak{o}(n)+\p$, where $\p$ is the linear space of symmetric operators in $\mathcal{L}$
and identify $\p^*$ with $\p$ via the Killing form. Then the moment map $\mm$ of the action of $G$ on $\mathbb{P}V$ can
be defined as in \eqref{eq:defmoment}, and one has the following result.

\begin{theorem}[{\cite{L1}, \cite{L3}}]\label{t:moment}
Let $\mu \in \mathcal{N} \setminus 0$ and let $\Ric$ be the Ricci endomorphism of the metric nilpotent Lie algebra
$(\mathcal{L}, \mu)$. Then

\emph{1.} $\mm([\mu]) = 4 \, \|\mu\|^{-2} \, \Ric$.

\emph{2.} Let $\n=(\Rn, \mu)$ be a nonabelian nilpotent Lie algebra. Choose an arbitrary inner product $\<\cdot,\cdot\>$
on $\Rn$. Then $\n$ is an Einstein nilradical if and only if
the function $F: \GL(n) \to \mathbb{R}$, the squared norm of the moment map, defined by
$F(g) = \|\mm(g.[\mu])\|^2 = 16 \, \|g.\mu\|^{-4} \, \Tr \, \Ric_{g.\mu}^2$
attains its minimum.
\end{theorem}

\section{Proof of Theorem~\ref{t:one}}
\label{s:proof}

The proof of Theorem~\ref{t:one} is a combination of Theorem~\ref{t:moment} and the results on convexity of
the image of the moment map of an orbit.

The first step in the proof is the observation that to check whether the condition of assertion 2
of Theorem~\ref{t:moment} is satisfied, one does not need the whole $\GL(n)$ orbit. Let $\n$ be a nilpotent Lie algebra
of dimension $n$. A derivation $\phi$ of $\n$ is called \emph{pre-Einstein}, if it is real (all the eigenvalues of
$\phi$ are real), semisimple, and
\begin{equation}\label{eq:pEtrace}
    \Tr (\phi \circ \psi) = \Tr \psi,  \quad \text{for any $\psi \in \Der(\n)$}.
\end{equation}
Note that the Einstein derivation, if it exists, satisfies \eqref{eq:pEtrace}, up to a nonzero multiple, as easily
follows from \eqref{eq:riccinil} and \eqref{eq:ricn}.
By \cite[Proposition 1]{N2}, a pre-Einstein derivation always exists, is unique up to
conjugation, and all its eigenvalues are rational. The main advantage of the pre-Einstein derivation in the study
of Einstein nilradicals lies the fact
that if $\n$ is an Einstein nilradical, then the Einstein derivation is a positive multiple of $\phi$ (up to
conjugation). Thus finding a pre-Einstein derivation for a given nilpotent Lie algebra $\n$, one immediately gets a
substantial portion of information on the nilsoliton inner product on $\n$, if the latter exists: for instance, in
our case, when all the eigenvalues of $\phi$ are simple, we already have an orthogonal (but not orthonormal!) basis
of Ricci eigenspaces of a (potentially existing) nilsoliton inner product. On the other hand, a
given algebra $\n$ is not an Einstein nilradical, if its pre-Einstein derivation fails to have all its eigenvalues
positive, or if the endomorphism $\ad_\phi$ of $\Der(\n)$ has nonpositive eigenvalues
(see the proof of Theorem 3 of \cite{N2}).

Let $Z(\phi) \subset \GL(n)$ be the centralizer of $\phi$ in $\GL(n)$, and $Z_0(\phi)$ its
identity component. %\textbf{(better: $Ad_G \phi = \phi$?)}
The group $Z_0(\phi)$ is isomorphic to $\prod_i \GL^+(d_i)$, where $d_1, \ldots, d_p$ are the multiplicities of
the eigenvalues of $\phi$ and $\GL^+(d) = \{M \in \GL(d) \, : \, \det M > 0\}$.

The following Lemma is essentially contained in \cite[Theorem 4.3]{L3} (see also \cite[Theorem 6.15, Lemma 6.14]{H}).
Note that at this stage we do not use the assumption that all the eigenvalues of $\phi$ are simple.

\begin{lemma} \label{l:centralizer}
Let $\n=(\Rn, \mu)$ be a nilpotent Lie algebra, with $\phi$ the pre-Einstein derivation. Let $Z_0(\phi) \subset \GL(n)$
be the identity component of the centralizer of $\phi$ in $\GL(n)$ and $\<\cdot,\cdot\>$ be an arbitrary inner product
on $\Rn$ with respect to which $\phi$ is symmetric. The algebra $\n$ is an Einstein nilradical if and only if
the function $F: Z_0(\phi) \to \mathbb{R}$ defined by $F(g) = \|\mm(g.[\mu])\|^2$ attains its minimum.
\end{lemma}

\begin{proof}
The claim follows from \cite[Theorem 4.3]{L3}, if we choose $G_\gamma$ to be $Z_0(\phi)$, and $\mathcal{C}$ to be the set
of inner products on $\n$ with respect to which $\phi$ is symmetric. Then $\Ric^\gamma$, the projection of $\Ric$ to the
Lie algebra $\z(\phi)$ of $Z_0(\phi)$, coincides with $\Ric$ by \cite[Lemma 2.2]{H}, as $\phi$ is a symmetric derivation.

The ``only if" part uses the fact that if a nilpotent Lie algebra is an Einstein nilradical, then its Einstein
derivation is proportional to a pre-Einstein derivation.
\end{proof}

Now suppose that all the eigenvalues $\la_i$ of $\phi$ are simple. Let $e_i$ be a basis of eigenvectors of $\phi$ and
let $[e_i, e_j] = \sum_{k=1}^n C_{ij}^k e_k$. The number $C_{ij}^k$ can be nonzero only if $\la_i+\la_j = \la_k$,
in particular, for every pair $(i, j)$, at most one of the $C_{ij}^k$ is nonzero.
The group $Z_0(\phi)$ is abelian and is isomorphic to $\GL^+(1)^n$ acting as follows: an element
$g=(e^{x_1}, \ldots, e^{x_n})$ sends $e_i$ to $e^{-x_i}e_i$. The corresponding action on $V$ is given by
$C_{ij}^k \to e^{x_i+x_j-x_k}C_{ij}^k$. Fix an inner product on $\n$ such that $\<e_i, e_j\> = \K_{ij}$.

The moment map $\mm$ acts to the space $\z^*(\phi)$ of diagonal matrices
with respect to the basis $\{e_i\}$. Identify $\z^*(\phi)$ with $\Rn$, with an inner product $(\cdot, \cdot)$ induced
by the Killing form of $\mathfrak{gl}(n)$ and with the orthonormal basis $\{f_i\}$ ($f_i$ corresponds to the matrix
having $1$ as its $(i,i)$-th entry and zero elsewhere).

Define $\mathbf{F}=\{\a_{ij}^k= f_i+f_j-f_k: C_{ij}^k \ne 0\} \subset \Rn$.
The claim of Theorem~\ref{t:one} immediately follows from Lemma~\ref{l:centralizer} and the following lemma:

\begin{lemma}\label{l:moment}
Let $\n=(\Rn, \mu)$ be a nilpotent Lie algebra whose pre-Einstein derivation $\phi$ has all its eigenvalues simple.
For an arbitrary inner product on $\n$ with respect to which $\phi$ is symmetric,
$$
\mm(Z_0(\phi).[\mu]) = (\Conv(\mathbf{F}))^0,
$$
the interior of the convex hull of $\mathbf{F}$.
\end{lemma}
\begin{proof}
By \cite[Proposition 4.1]{HS}, the set $\mm(Z_0(\phi).[\mu])$ is an open convex subset of an affine subspace of $\Rn$.
To prove the lemma it therefore suffices to show that $\mm(Z_0(\phi).[\mu]) \subset \Conv(\mathbf{F})$ and that
$\overline{\mm(Z_0(\phi).[\mu])} \supset \mathbf{F}$.

Let $X = (x_1, \ldots, x_n) \in \z(\phi) = \Rn$ and $g= \exp X = (e^{x_1}, \ldots, e^{x_n}) \in Z_0(\phi)$. Then
\begin{equation}\label{eq:moment}
\mm(g.[\mu])=\frac{\sum_{\a_{ij}^k \in \mathbf{F}} e^{2(\a_{ij}^k,X)} (C_{ij}^k)^2\a_{ij}^k}
{\sum_{\a_{ij}^k \in \mathbf{F}} e^{2(\a_{ij}^k,X)}(C_{ij}^k)^2}.
\end{equation}
It follows that $\mm(Z_0(\phi).[\mu]) \subset \Conv(\mathbf{F})$, as for every
$g \in Z_0(\phi), \quad \mm(g.[\mu])$ is the center of mass of the set $\mathbf{F}$, with the positive masses
$e^{2(\a_{ij}^k,X)}(C_{ij}^k)^2$ placed at the vertices $\a_{ij}^k$.
Moreover, $\overline{\mm(Z_0(\phi).[\mu])}$ contains all the vertices of $\mathbf{F}$. Indeed, let
$\a_{ij}^k \in \mathbf{F}$ and let $X = f_i+f_j$. Then $(\a_{ij}^k, X)=2$ and $(\a_{ls}^r, X) < 2$ for any
other vertex $\a_{ls}^r \in \mathbf{F}$. By \eqref{eq:moment},
$\lim_{t\to \infty} \mathbf{m}(\exp(tX). [\mu]) = \a_{ij}^k$.
\end{proof}

\begin{remark} \label{rem:multiplicity}
A direct generalization of Theorem~\ref{t:one} to the case when the pre-Einstein derivation $\phi$ has eigenvalues of
higher multiplicities works only as a necessary condition. If $\la_1, \ldots, \la_p$ are the eigenvalues of $\phi$,
with $\n_1, \ldots, \n_p$ the corresponding eigenspaces, then for every pair $(i,j), \quad [\n_i,\n_j]$ is either zero,
or lies in some eigenspace $\n_k$. Defining $\mathbf{F}$ as the subset of $\mathbb{R}^p$ consisting of the
vectors $\a_{ij}^k=f_i+f_j-f_k$ such that $[\n_i, \n_j] \subset \n_k,\; [\n_i, \n_j] \ne 0$ we get a necessary condition
for $\n$ to be an Einstein nilradical similar to that of \cite[Lemma 1]{N1}.

The reason why the ``if" part of Theorem~\ref{t:one} fails in this case is because Lemma~\ref{l:moment} is no
longer true. In general,
let $G$ be a reductive Lie group acting on a linear space $V$ and let $\g= \k \oplus \p$
be the Cartan decomposition of its Lie algebra, with $\k$ the Lie algebra of a maximal compact subgroup
$K \subset G$. Let $\ag \in \p$ be a maximal subalgebra (which is always abelian, as $[\p,\p] \subset \k$).
The image of the moment map $\mathbf{m}$ of the action of $G$ on the projective space $\mathbb{P}V$ lies in $\p$.
One has two sorts of the general convexity results for the image of the $G$-orbit of a point $x \in \mathbb{P}V$
under the moment map $\mathbf{m}$: the projection of $\mathbf{m}(G.x)$ to $\ag$ is convex (the Kostant Theorem)
and the intersection of $\mathbf{m}(\overline{G.x})$ with the positive Weyl chamber $\ag_+ \subset \ag$ is
convex \cite[Corollary 7.1]{S}.

If all the eigenvalues of $\phi$ are simple, the group $G= \prod \GL^+(1) \cong \Rn$ is abelian,
so $\k =0, \g=\p=\ag$ and the Weyl group is trivial,
hence $\ag_+ = \ag (= \Rn)$, which implies that the image of the orbit is convex in $\ag$. In general, however, this is
not true even in very simple cases, as the following example shows
(some examples of that sort in the settings different from ours can be found in \cite[Chapter 8]{S}).
\end{remark}

\begin{example}
Consider a $(2p+1)$-dimensional two-step nilpotent Lie algebra $\n$ given by the relations
$[e_1, e_i]=e_{i+p}, \; i=2, \ldots p+1$ (note that such an algebra is an Einstein nilradical
by \cite[Theorem 4.2]{L2}). A pre-Einstein derivation $\phi$ can be taken as
$\phi(e_1) = \frac{2}{p+2} e_1,\; \phi(e_i) = \frac{p+1}{p+2} e_i, \, i=2,\ldots, p+1,\;
\phi(e_j) = \frac{p+3}{p+2} e_j, \, j=p+2,\ldots, 2p+1$
(in fact, $\frac{2}{p+2} \phi$ is an Einstein derivation, if we choose an inner product on $\n$ in such a way that
the vectors $e_i$ are orthonormal). The component of the identity of the centralizer of $\phi$
is the group $G=\GL^+(1) \times \GL^+(p) \times \GL^+(p)$.
Introduce the inner product on $\mathbb{R}^{2p+1}$, the linear space of $\n$ by requiring that the basis $e_i$ is
orthonormal. Denote $\mu$ the bracket defining $\n$. Then for $g = (t, g_1, g_2) \in G$ the bracket $g.\mu$ is given
by $g.\mu(e_1,e_i) = t^{-1} g_2 [e_1, g_1^{-1} e_i], \; 2 \le i \le p+1$ (and $g.\mu(e_i,e_j) = 0$ for all the other
pairs with $i < j$). For any $g \in G$ there exist an $h \in \GL^+(p)$ (acting on the space
$\Span (e_{p+2}, \ldots, e_{2p+1})$) such that $g.\mu(e_1,e_i) = h e_{i+p}$ for $2 \le i \le p+1$. Such an $h$ is
uniquely determined by $g \in G$ and the map from $G$ to $\GL^+(p)$ sending $g$ to $h$ is onto.

The moment map for the action of $G$ on $\mu$ is given by
$$
\mathbf{m}(g.[\mu])=
\begin{bmatrix}
  -2 & 0 & 0 \\
  0 & -2 h^th \cdot (\Tr h^t h )^{-1}& 0 \\
  0 & 0 & 2 hh^t \cdot (\Tr hh^t)^{-1} \\
\end{bmatrix}
$$
The intersection of $\mathbf{m}(g.[\mu])$ with $\ag$, the set of diagonal matrices, is the set
$\{\diag(-2,-\la_1, \ldots, -\la_p,$
$\la_{\sigma(1)}, \ldots, \la_{\sigma(p)}) \, : \, \la_i > 0, \, \sum \la_i = 2, \, \sigma \in S_p\}$,
where $S_p$ is the symmetric group of order $p$. So $\mathbf{m}(g.[\mu]) \cap \ag$ is the union of
$p! \; (p-1)$-dimensional open simplices and is not convex.
For instance, for $p=2, \quad \mathbf{m}(g.[\mu]) \cap \ag$ is the union of two diagonals of a square.
It is easy to see, however, that $\mathbf{m}(g.[\mu]) \cap \ag_+$ is convex (and is a simplex minus some faces).

\end{example}

\bigskip

The class of nilpotent Lie algebras whose pre-Einstein derivation has all its eigenvalues simple is quite large.
One example, the graded filiform algebras, will be considered in details
in the next section. However, there are many other algebras with such a property. For instance, a pre-Einstein
derivation of a nilpotent algebra with codimension one abelian ideal $\ag$ (see \cite[Section 4]{L2}) has all its
eigenvalues simple, provided the operator $\ad_X$ (where $X \notin \ag$) has two Jordan blocks of
dimensions $p_1$ and $p_2$, with $p_1+p_2$ an odd number. Another example is a two-step nilpotent algebra given by
the relations
$[e_1, e_2] = e_7, \, [e_1, e_3] = e_8, \,[e_2, e_3] = [e_4, e_5] = e_9$, $[e_3, e_4] = [e_1, e_6] = e_{10}$.
We give yet another example below.
\begin{example}
Consider a family of eight-dimensional nilpotent Lie algebras given with respect to a basis $e_i,\; i=1,\ldots 8$,
by the relations $[e_i,e_j] = c_{ij} e_{i+j+1}$, where all the $c_{ij}$'s with $i<j,\, i+j \le 7$, are nonzero.
Any such algebra is isomorphic to exactly one of the algebras $\n_t(8), \; t \ne 0, 1$, defined by
$[e_1, e_i] = e_{i+2}, \; i = 2, \ldots, 6, \;
[e_2, e_3] = e_6,\, [e_2, e_4] = e_7,\, [e_2, e_5] = t e_8,\, [e_3, e_4] = (t-1) e_8$. The pre-Einstein derivation
for each of the $\n_t(8), \; t \ne 0, 1$, is positively proportional to the derivation sending
every $e_i$ to $(i+1)e_i$, hence has all its eigenvalues simple. The routine check shows that the set $S$ of positive
solutions $v$ of \eqref{eq:ut=1} (see Remark~\ref{rem:payne}) is nonempty. In fact, all the $\n_t(8), \; t \ne 0, 1$,
share the same $S$, which is the interior of a triangle in $\mathbb{R}^9$. By Theorem~\ref{t:one}, each of the algebras
$\n_t(8), \; t \ne 0, 1$ is an Einstein nilradical.
\end{example}

\section{Filiform $\mathbb{N}$-graded algebras and Einstein nilradicals}
\label{s:fili}

In this section, we apply Theorem~\ref{t:one} to study filiform algebras.
A \emph{filiform} algebra is a nilpotent algebra for which the descending central series
$\n_0=\n,\; \n_{i+1} = [\n, \n_i],\, i \ge 0$, has the
maximal possible length for the given dimension, namely $\n_{n-1} \ne 0$, where $n = \dim \n$.

We are primarily interested in those filiform algebras which are Einstein nilradicals. Every such algebra must admit
an $\mathbb{N}$-gradation.

Any filiform algebra $\n$ is generated by two elements, so its rank is at most two.
Most of filiform algebras of dimension $n \ge 8$ are characteristically nilpotent \cite{GH1}, that is, have
rank zero. Those algebras do not admit any gradation at all, hence cannot be Einstein nilradicals
(cf. \cite[Section 9.1]{P}); in fact, they cannot even be the nilradicals of anything other than themselves.

The case $\rk \, \n = 2$ is completely settled by the following two facts. First of all, by the result of \cite{V, GH1},
there are only two filiform Lie algebras of rank two, namely
\begin{align}\label{eq:m0}
&\m_0(n) \; &[e_1, e_i] = e_{i+1}, \quad &i=2, \ldots n-1, \\
&\m_1(n) \; &[e_1, e_i] = e_{i+1}, \quad &i=2, \ldots n-2, \quad [e_i,e_{n-i+1}] = (-1)^i e_n, \; i = 2, \ldots , n - 1,
\label{eq:m1}
\end{align}
where the dimension $n$ of $\m_1(n)$ must be even. Secondly, both of these algebras are Einstein nilradicals. For
$\m_0(n)$, this is proved in \cite[Theorem 4.2]{L2} (see also \cite[Theorem 27]{P}), for $\m_1(n)$ in
\cite[Theorem 37]{P}. Note that the pre-Einstein derivation for both $\m_0$ and $\m_1$ is simple, so the fact
that they are Einstein nilradicals can be also deduced from Theorem~\ref{t:one}.

Less is known in the case $\rk \, \n = 1$. As it follows from \cite[Th\'{e}or\`{e}me 2]{GH1} (see also
\cite[Section~3.1]{GH2}), there are two series of classes of rank one filiform algebras admitting a positive gradation:
\begin{itemize}
  \item the $n$-dimensional algebras of the class $A_r,\; 2\le r \le n-3$ are given by the relations
  \begin{equation*}
    [e_1, e_i] = e_{i+1}, \quad i=2, \ldots n-1, \quad [e_i, e_j] = c_{ij} e_{i+j+r-2} , \quad i,j \ge 2, i+j \le n+2-r,
  \end{equation*}
  with the gradation $1, r, r+1, \ldots, n+r-2$ (the corresponding derivation $\phi$ is defined by
  $\phi(e_1) = e_1$, $\phi(e_i) = (i+r-2) e_i,\;i \ge 2$).
  \item the $n$-dimensional algebras of the class $B_r,\; 2\le r \le n-3$, $n$ is even, are given by the relations
  \begin{equation}\label{eq:Brn}
    [e_1, e_i] = e_{i+1}, \quad i=2, \ldots n-2, \quad [e_i, e_j] = c_{ij} e_{i+j+r-2} , \quad i,j \ge 2, i+j \le n+2-r,
  \end{equation}
  with the gradation $1, r, r+1, \ldots, n+r-3, n+2r-3$ (the corresponding derivation $\phi$ is defined by
  $\phi(e_1) = e_1,\; \phi(e_i) = (i+r-2) e_i,\; 2 \le i \le n-1,\; \phi(e_n) = (n+2r-3) e_n$).
\end{itemize}
In order to get an algebra of rank precisely one, one requires that not all of the $c_{ij}$'s above are zeros.
The main difficulty in classifying the algebras from $A_r$ and $B_r$ lies in the fact that the
numbers $c_{ij}$ must satisfy the Jacobi equation. Note that for any $n \ge 5$ and any $2\le r \le n-3$,
there exists, for instance, an algebra from $A_r$ with all the $c_{ij}$'s nonzero. The complete classification is
known only for the algebras of the class $A_2$ (\cite{M,CJ}, the complex case was earlier done in \cite{AG}).

Note that although the ground field in \cite{GH1, GH2} is $\mathbb{C}$, the examination of the proof shows that the
classification of filiform algebras of rank one works for $\mathbb{R}$ without any changes (actually, for any field
of infinite characteristics; note, however, that the algebra $\b_{\pm}(12)$ from Table~\ref{tableb2n} below is not
defined over $\mathbb{Q}$).

Every algebra $\n$ from $A_r$ or $B_r$ has only one semisimple derivation, which is automatically a
pre-Einstein derivation (up to conjugation and scaling). As all the eigenvalues of it are simple
(they are proportional to $(1, r, r+1, \ldots, n+r-2)$ for $A_r$ and to $(1, r, r+1, \ldots, n+r-3, n+2r-3)$ for $B_r$),
the question of whether or not $\n$ is an Einstein nilradical is answered by Theorem~\ref{t:one}.

Note that the affine space $L$ in Theorem~\ref{t:one} (and hence the projection $p$ of the origin of $\Rn$ to it) is
the same for all the $n$-dimensional algebras of each of the classes $A_r$ and $B_r$ (although the set $\mathbf{F}$
depends on the particular algebra). The explicit form of $p$ can be easily found, see e.g. \eqref{eq:pA2n} for $A_2$.

With the classification of algebras of the class $A_2$ in hands, we classify algebras of class $B_2$
(Table~\ref{tableb2n} in Section~\ref{ss:b2n}); it appears that there are only six of them. Then we apply
Theorem~\ref{t:one} to find all the Einstein nilradicals in the classes $A_2$ and $B_2$, which proves
Theorem~\ref{t:fili}.

\subsection{Algebras of the class $\mathbf{A_2}$}
\label{ss:a2n}

According to the classification given in \cite[Theorem 5.17]{M}, the class $A_2$ consists of five infinite
series and five one-parameter families in dimensions $7 \le n \le 11$. More precisely, there are two infinite
series, $\m_2(n)$ and $\mathcal{V}_n$, defined for all $n$, two others, $\m_{0,1}(n)$ and $\m_{0,3}(n)$, defined
for odd $n$, and one, $\m_{0,2}(n)$, defined for even $n$. The tables below give the commuting relations for the
algebras from $A_2$ (they slightly differ from the ones from \cite{M}: first of all, we remove the algebra $\m_0$,
as it is of rank two; secondly, we change the lower bounds for the dimensions, so in our tables some lower-dimensional
algebras from different families could be isomorphic).

\begin{table}[h]
\setlength{\extrarowheight}{2pt}
\begin{center}
\begin{tabular}{|m{3cm}|m{11cm}|}
\hline
$\m_2(n), \, n \ge 5$ & $[e_1, e_i] = e_{i+1},$ \hfill $i = 2, \ldots , n-1$ \newline
$[e_2, e_i] = e_{i+2},$ \hfill $i = 3, \ldots , n-2$  \\
\hline
$\mathcal{V}(n), \, n \ge 4$ & $[e_i, e_j ] = (j-i) e_{i+j},$ \hfill $i + j \le n$ \\
\hline
$\m_{0,1}(n), $ \newline $n=2m+1,\, n \ge 7$ & $[e_1, e_i] = e_{i+1},$ \hfill $i = 2, \ldots , n-1$ \newline
$[e_l, e_{n-l}] = (-1)^{l+1} e_n,$ \hfill $l = 2, \ldots , m$ \\
\hline
$\m_{0,2}(n),$ \newline $n=2m+2,\, n \ge 8$ & $[e_1, e_i] = e_{i+1},$ \hfill $i=2, \ldots , n-1$ \newline
$[e_l, e_{n-1-l}] = (-1)^{l+1} e_{n-1},$ \hfill $l=2, \ldots , m$ \newline
$[e_j , e_{n-j}] = (-1)^{j+1}(m-j+1)e_n,$ \hfill $j=2, \ldots , m$ \\
\hline
$\m_{0,3}(n),$ \newline $n=2m+3,\, n \ge 9$ &
$[e_1, e_i] = e_{i+1},$ \hfill $i=2, \ldots , n-1$ \newline
$[e_l, e_{n-2-l}] = (-1)^{l+1}e_{n-2},$ \hfill $l=2, \ldots , m$ \newline
$[e_j , e_{n-1-j}] = (-1)^{j+1}(m-j+1)e_{n-1},$ \hfill $j=2, \ldots , m$ \newline
$[e_k, e_{n-k}] = (-1)^k ((k-2)m - (k-2)(k-1)/2)e_n,$ \hfill $k=3, \ldots , m+1$ \\
\hline
\end{tabular}
\caption{Infinite series of algebras of the class $A_2$.}\label{tablea2}
\bigskip
\begin{tabular}{|m{1cm}|m{12.9cm}|}
\hline
$\g_\a(7)$ & $[e_1, e_i] = e_{i+1}$, \quad $i = 2, \ldots , 6$ \newline
$[e_2, e_3] = (2+\a) e_5,\, [e_2, e_4] = (2+\a) e_6,\, [e_2, e_5] = (1+\a) e_7,\, [e_3, e_4] = e_7$\\
\hline
$\g_\a(8)$ & relations for $\g_\a(7)$ and \newline
$[e_1, e_7 ] = e_8, \, [e_2, e_6 ] = \a e_8, \, [e_3, e_5 ] = e_8$\\
\hline
$\g_\a(9)$ & relations for $\g_\a(8)$ and \newline
$[e_1, e_8] = e_9, \, [e_2, e_7] = \frac{2 \a^2 + 3 \a - 2}{2\a+5} e_9, \,
[e_3, e_6] = \frac{2 \a + 2}{2\a+5} e_9, \, [e_4, e_5] = \frac{3}{2\a+5} e_9$, \hfill $\a \ne -\frac52$\\
\hline
$\g_\a(10)$ & relations for $\g_\a(9)$ and \newline
$[e_1, e_9] = e_{10},\, [e_2, e_8] = \frac{2 \a^2 + \a - 1}{2\a+5} e_{10}, \,
[e_3, e_7] = \frac{2 \a - 1}{2\a+5} e_{10}, \, [e_4, e_6] = \frac{3}{2\a+5} e_{10}$, \hfill $\a \ne -\frac52$ \\
\hline
$\g_\a(11)$ & relations for $\g_\a(10)$ and \newline
$[e_1, e_{10}] = e_{11},\, [e_2, e_9] = \frac{2 \a^3 + 2 \a^2 + 3}{2(\a^2+4\a+3)} e_{11}, \,
[e_3, e_8] = \frac{4 \a^3 + 8 \a^2 - 8 \a - 21}{2(\a^2+4\a+3)(2\a+5)} e_{11}$, \newline
$[e_4, e_7] = \frac{3(2 \a^2 + 4 \a + 5)}{2(\a^2+4\a+3)(2\a+5)} e_{11}, \,
[e_5, e_6] = \frac{3(4 \a + 1)}{2(\a^2+4\a+3)(2\a+5)} e_{11}$,  \hfill $\a \ne -3, -\frac52, -1$\\
\hline
\end{tabular}
\caption{One-parameter families of algebras of the class $A_2$.}\label{tablega}
\end{center}
\end{table}

As we are interested only in determining whether a given algebra is an Einstein nilradical, by
Theorem~\ref{t:one}, we need only the set $\mathbf{F}$ for each of the algebras, not the actual structural
coefficients.
The vector $p$, the projection of the origin of $\Rn$ to $L$, is given by
\begin{equation}\label{eq:pA2n}
    p_i= \frac{2}{n(n-1)}(2n+1-3i), \quad i=1, \ldots , n.
\end{equation}

The proof goes on the case-by-case basis.

First of all, neither of the algebras
$\m_2(n), \; \m_{0,1}(n), \, n=2m+1, \; \m_{0,2}(n), \, n=2m+2, \; \m_{0,3}(n), \; n=2m+3$ is an Einstein nilradical,
when $n \ge 8$. The easiest way to see that is to produce a vector $a \in \Rn$ such that $(a, \a_{ij}^k) \ge 0$, for
all $\a_{ij}^k \in \mathbf{F}$, but $(a, p) < 0$ (this implies that $p \notin \Conv(\mathbf{F}))$. Such a vector $a$
can be taken as $a_1=n-2, \, a_2=2(n-2), \, a_i =i(n-2)-n(n+1)/2, \, 3 \le i \le n$, for $\m_2(n)$, as
$(1,1-m,2-m, \ldots, m-2, m-1,-1)^t$ for $\m_{0,1}(n)$, as $a=(1, 1-m, 2-m, \ldots, m-2, m-1,1,0)^t$ for
$\m_{0,2}(n), \, n=2m+2$, and as
$a_1=n+2, \, a_{n-2}=-n-4, \,  a_{n-1}=-2, \,  a_n=n, \, a_i =i(n+2)-\frac{n(n+1)}{2}$, $3 \le i \le n-3$,
for $\m_{0,3}(n)$.

The only remaining algebra from Table~\ref{tablea2}, the algebra $\mathcal{V}(n)$, is an Einstein nilradical for any
$n \ge 3$. This follows from the fact that
$p=\frac{2}{n(n-1)} (2 \sum_{1\le i < j, \, i+j \le n} \a_{ij}^{i+j}+\a_{1m}^{m+1}+\sum_{i=1}^{m-1} \a_{i,i+2}^{2i+2})$,
where $m =[n/2]$.

Apart for some finite number of exceptional values of $\a$, for every $n=7, \ldots, 11$, the set $\mathbf{F}$ for the
algebras $\g_\a(n)$ from Table~\ref{tablega} is the same as for the corresponding algebra $\mathcal{V}(n)$, so each
of them is an Einstein nilradical. We treat the exceptional values below by either giving the coefficient vector
$c=(c_{ij}^k)$ of a convex linear combination $p=\sum c_{ij}^k \a_{ij}^k$ (the vectors $\a_{ij}^k$ are always ordered
lexicographically), or otherwise, by showing that some coefficient of any such linear combination is nonpositive.
According to Theorem~\ref{t:one}, the corresponding algebra is an Einstein nilradical in the former case, and is not
in the latter one.

The exceptional values for $\g_{\a}(8)$ are $\a= -2, -1, 0$. The algebra $\g_{-2}(8)$ is not an Einstein nilradical, as
for any linear combination of the vectors from $\mathbf{F}$ representing $p, \; c_{16}^7 < 0$.
Both $\g_{0}(8)$ and $\g_{-1}(8)$ are Einstein nilradicals, the coefficient vector $c$ can be taken as
$\frac{1}{28} (4, 2, 3, 1, 2, 2, 3, 2, 2, 2, 5)^t$ for $\a=0$ and as
$\frac{1}{28}(2, 2, 2, 4, 3, 1, 3, 3, 2, 3, 3)^t$ for $\a=-1$.

The exceptional values for $\g_{\a}(9)$ are $\a= -2, -1, 0$. The algebra $\g_{-2}(9)$ is not an Einstein nilradical, as
for any linear combination of the vectors from $\mathbf{F}$ representing $p, \; c_{16}^7 < 0$.
All three algebras $\g_{-1}(9)$, $\g_{0}(9)$, and $\g_{1/2}(9)$ are Einstein nilradicals, the coefficient vector $c$ can
be taken as $\frac{1}{36} (2, 2, 1, 2, 2, 2, 5, 3, 3, 4, 1, 3, 4, 2)^t$, as
$\frac{1}{72} (6, 6, 5, 3, 5, 6, 1, 5, 5, 5, 5, 5, 5, 5, 5)^t$, and as
$\frac{1}{36} (3, 1, 2, 4, 1, 3, 2$, $4, 2, 2, 2, 2, 2, 4, 2)^t$, respectively.

The exceptional values for $\g_{\a}(10)$ are $\a= -2, -1, 0, \frac12$. Neither of the algebras
$\g_{-2}(10)$, $\g_{-1}(10)$, and $\g_{1/2}(10)$ is an Einstein nilradical, as any linear combination of the vectors
from $\mathbf{F}$ representing $p$ has some of the coefficient nonpositive (specifically, for
$\g_{-2}(10)$, we have $c_{16}^7 + c_{19}^{10} + c_{46}^{10} = 0$, for $\g_{-1}(10)$, $c_{14}^5 + c_{17}^8 = 0$,
for $\g_{1/2}(10)$, $c_{14}^5 + c_{16}^7 + c_{34}^7 = 0$). The algebra $\g_{0}(10)$ is an Einstein nilradical,
the coefficient vector $c$ can be taken as $\frac{1}{45} (4, 3, 2, 1, 2, 2, 2, 2, 2, 3, 2, 2, 2, 2, 5, 2, 2, 2, 3)^t$.

The exceptional values for $\g_{\a}(11)$ are $\a= -2, -1/4, 0, 1/2, \a_1, \a_2$, where $\a_1$
is the unique real root of $2 \a^3 + 2 \a^2 + 3$, and $\a_2$ is the unique real root of $4 \a^3 + 8 \a^2 - 8 \a - 21$.
All these algebras are Einstein nilradicals: the coefficient vectors $c$ can be taken as
%\newline
$\frac{1}{220} (22, 12, 15, 12, 1, 8, 8, 1, 1, 8, 8, 8, 22, 14, 12, 17, 15, 8$, $14, 5, 8, 1)^t$,
$\frac{1}{55} (2, 3, 2, 3, 2, 2, 2, 2, 2, 2, 1, 4, 2, 2, 2, 2, 2, 3, 2, 2, 2, 2, 3, 4)^t$,
$\frac{1}{55} (2, 2, 2, 2, 2, 4, 1, 3, 2, 2, 1, 4, 4$, $2, 2, 4, 2, 2, 2, 2, 2, 2, 2, 2)^t$,
$\frac{1}{110} (9, 10, 1, 6, 1, 4, 4, 4, 1, 7, 6, 4, 4, 4, 1, 6, 5, 8, 7, 11, 6, 1)^t$,
$\frac{1}{55} (3, 3, 2, 3, 2, 2, 1$, $2, 2, 3, 3, 2, 2, 2, 2, 2, 3, 1, 2, 3, 2, 3, 2, 3)^t$,
and $\frac{1}{55} (1, 4, 2, 3, 2, 1, 3, 2, 2, 2, 3, 2, 2, 2, 2, 3, 2, 3, 2, 2, 2, 3, 3, 2)^t$, \linebreak
respectively.

\subsection{Algebras of the class $\mathbf{B_2}$}
\label{ss:b2n}

In this section, we classify the filiform algebras of the class $B_2$ and show that all of them are Einstein nilradicals,
hence proving assertion 2 of Theorem~\ref{t:fili}.

An $n$-dimensional algebra $\n \in B_2$ is defined by the relations \eqref{eq:Brn}, with $r=2$. Any such algebra
admits a derivation with the eigenvalues $1,2,\ldots, n-1, n+1$ (the corresponding eigenvectors are the $e_i$'s).
We require that $\rk \, \n = 1$, that is, at least one of the $c_{ij}$ is nonzero. This sorts out
the rank two algebra $\m_1(n)$ given by \eqref{eq:m1}.

We first prove the classification part of assertion 2 of Theorem~\ref{t:fili}. Namely, we show that $B_2$ consists
of the six algebras given in Table~\ref{tableb2n}. Each of those algebras has an even dimension $n = 2m$ and is a
cental extension of some $n-1$-dimensional algebra of the class $A_2$ by the cocycle
$$
\omega = a_2 e_2^* \wedge e_{n-1}^* + a_3 e_3^* \wedge e_{n-2}^* + \ldots + a_m e_m^* \wedge e_{m+1}^*.
$$

\begin{table}[h]
\setlength{\extrarowheight}{2pt}
\begin{center}
\begin{tabular}{|m{1.25cm}|m{3.2cm}|m{9.7cm}|}
\hline
Algebra & Extension of & {\centering Relations}\\
\hline
$\b(6)$ & $\m_2(5)$ & relations for $\m_2(5)$ and
$[e_i, e_{7-i}] = (-1)^i e_6, \hfill i=2,3$\\
\hline
$\b(8)$ & $\g_{-5/2}(7)$ &
relations for $\g_{-5/2}(7)$ and
$[e_i, e_{9-i}] = (-1)^i e_6, \hfill i=2,3,4$\\
\hline
$\b_1(10)$ & $\g_{-1}(9)$ &
relations for $\g_{-1}(9)$ and
$[e_i, e_{11-i}] = (-1)^i e_{10}, \hfill i=2,3,4,5$\\
\hline
$\b_2(10)$ & $\g_{-3}(9)$ &
relations for $\g_{-3}(9)$ and
$[e_i, e_{11-i}] = (-1)^i e_{10}, \hfill i=2,3,4,5$\\
\hline
$\b_{\pm}(12)$ & $\g_{\a}(11), \; \a=\frac{-4 \pm \sqrt10}{2}$ &
relations for $\g_{\frac{-4 \pm \sqrt10}{2}}(12)$ and %\newline
$[e_i, e_{13-i}] = (-1)^i e_{12}, \hfill i=2,\ldots,6$\\
\hline
\end{tabular}
\caption{Algebras of the class $B_2$.}\label{tableb2n}
\end{center}
\end{table}

\begin{proof}[Proof of the classification]
Let $\n$ be a filiform algebra of dimension $n$ from the class $B_2$. Then $\n$ is of rank one and admits a derivation
with the eigenvalues $1,2,\ldots, n-1, n+1$. Let $\{e_i\}$ be the corresponding basis of the eigenvectors. Then
$\z(\n)$, the center of $\n$, is $\mathbb{R} e_n$, and the quotient algebra $\n'=\n/\z(\n)$ is a filiform
$n-1$-dimensional algebra admitting
a derivation with the eigenvalues $1,2, \ldots, n-1$. The corresponding eigenvectors are the images of the
vectors $e_i, \, i=1, \ldots, n-1$, under the natural projection $\pi:\n \to \n'$. With a slight abuse of notation,
we will still denote them $e_i$.

The algebra $\n'$ is either isomorphic to $\m_0$, or is one of the algebras from $A_2$. Up to scaling,
we can assume that $[e_1, e_i]=e_{i+1}$, for all $i=2, \ldots, n-2$,
and $[e_i, e_j] = c_{ij} e_{i+j}$ when $i+j \le n-1$, or zero otherwise.

Given $\n'$, one can construct $\n$ as a central extension of $\n'$, so that $\n = \n' \oplus \mathbb{R} e_n$
(as a linear space), with the Lie brackets given by $[\n, e_n] = 0, \; [X,Y]=[X,Y]_1 + \omega(X, Y) e_n$,
for $X, Y \in \n'$, where $[X,Y]_1$ is the bracket of $X$ and $Y$ in $\n'$ and $\omega \in \Lambda^2(\n')$.
The Jacobi equations are equivalent to the fact that $\omega$ is a $2$-cocycle, that is
$\sigma_{XYZ}(\omega([X,Y],Z)) = 0$, for any $X, Y, Z \in \n'$, where $\sigma$ is the sum of the cyclic permutations.
The fact that $\n$ is filiform implies that $\omega \ne 0$. Clearly, the proportional $\omega$'s yield isomorphic
algebras.

In order for $\n$ to admit the gradation $1,2, \ldots, n-1, n+1$, the cocycle $\omega$ must be of a very special form,
namely $\omega(e_i, e_j) = 0$ for all $i,j =1, \ldots, n-1$, unless $i+j = n+1$. It follows that
$\omega= \sum_{i=2}^{n-1} a_i e_i^*\wedge e_{n+1-i}^*$, with $a_{n+1-i} = -a_i$. The cocycle condition with
$X=e_1, \, Y=e_i, \, Z=e_{n-i}$, $2 \le i \le n-2$, implies $a_{i+1} = -a_i$. As $\omega \ne 0$, $n$ must be even,
and up to scaling, we can take $\omega = \sum_{i=2}^m (-1)^i e_i^*\wedge e_{n+1-i}^*$, where $2m = n$.

For a triple $X=e_i, \, Y=e_j, \, Z= e_k$, the cocycle condition is nontrivial, only when $i+j+k = n+1$ and $i, j, k$ are
pairwise distinct. In such a case we get
\begin{equation}\label{eq:cocycle}
c_{ij} (-1)^k + c_{jk} (-1)^i + c_{ki} (-1)^j = 0, \quad \text{for all $i, j, k \ge 2$, with $i+j+k= n+1$.}
\end{equation}
Now, if $\n' = \m_0$, the condition \eqref{eq:cocycle} is clearly satisfied, the resulting algebra $\n$ is
isomorphic to $\m_1$ given by \eqref{eq:m1}. If not, then $\n'  \in A_2$, so it is one of the algebras from
Table~\ref{tablea2} or from Table~\ref{tablega}.

The direct check shows that the only odd-dimensional algebra from Table~\ref{tablea2} satisfying \eqref{eq:cocycle}
is the algebra $\m_2(5) \cong \mathcal{V}(5)$. It extends to $\b(6)$.

From among the algebras in Table~\ref{tablega}, only $\g_{-5/2}(7), \; \g_{-1}(9), \; \g_{-3}(9)$ and
$\g_{\a}(11),\; \a=\frac{-4 \pm \sqrt10}{2}$ satisfy \eqref{eq:cocycle}. The corresponding algebras of the
class $B_2$ are given in Table~\ref{tableb2n}.
\end{proof}

All the algebras of the class $B_2$ are Einstein nilradicals. It suffices to produce for each of them a
vector $v$ with positive coordinates satisfying \eqref{eq:ut=1}. Such a $v$ can be taken as
$\frac{1}{52} (13, 16, 12, 4, 13, 12)^t$,
$\frac{1}{221}(44, 11, 33, 79, 17, 48, 17, 21, 17, 17, 78, 17)^t, %\newline
\frac{1}{29}(2, 5, 4, 4, 4, 2, 4, 2, 4, 5, 4, 4, 4, 3, 4, 4, 3, 4)^t, \frac{1}{29}(3, 5, 3, 3$, \linebreak
$3, 3, 5, 3, 3, 3, 3, 3, 5, 4, 1, 1, 3, 2, 6, 4)^t$, %\newline
$\frac{1}{675} (33, 4, 93, 151, 105, 137, 8, 45, 20, 130,
45, 45, 45, 45, 45, 84, 45$, \linebreak
$45, 112, 45, 45, 45, 45, 45, 45, 45, 45, 45, 172, 45)^t$,
for the algebras $\b_6, \b_8, \b_1(10)$, $\b_2(10), \b_{\pm}(12)$, respectively.

%\textbf{
%0.check, no coeff in front in eq:moment? (shouldn't be)\newline
%5. all refs on def of moment contain a factor $i$ which we omit (to the first line of the proof of l2)\newline
%9. add ref to stokel?\newline
%10. men'she citirovanij l
%}

\end{document}